\documentclass{elsart}
\usepackage{mathrsfs}
\usepackage{epsfig,latexsym,amsfonts,amsmath,amssymb}

\begin{document}
\begin{frontmatter}

\newtheorem{theorem}[subsection]{Theorem}
\newtheorem{lemma}[subsection]{Lemma}
\newtheorem{conjecture}[subsection]{Conjecture}
\newtheorem{proposition}[subsection]{Proposition}
\newtheorem{definition}[subsection]{Definition}
\newtheorem{corollary}[subsection]{Corollary}
\newtheorem{example}[subsection]{Example}
\newtheorem{remark}[subsection]{Remark}
\newtheorem{de}[subsection]{Definition}

\renewcommand{\theequation}{\arabic{section}. \arabic{equation}}
\renewcommand{\thesection}{\arabic{section}}
\renewcommand{\thethm}{\arabic{section}.\arabic{thm}}
\renewcommand{\thefootnote}{\fnsymbol{footnote}}

\newcommand{\df}{\stackrel{\mbox{\rm def }}{=}}
\newcommand{\di}{\displaystyle}
\newcommand{\bl}[2]{{\left\langle #1 \:\:  \vrule \:\:  #2 \right\rangle}}
\newcommand{\vanish}[1]{}

\title{Cross-intersecting families and primitivity of symmetric systems}

\author[wj]{Jun Wang}\ead{jwang@shnu.edu.cn}\ and
\author[wj,zh]{Huajun  Zhang\thanksref{Y}}\ead{huajunzhang@zjnu.cn}
\address[wj]{ Department of Mathematics, Shanghai Normal University,
Shanghai 200234, P.R. China}
\address[zh]{Department of Mathematics,
 Zhejiang Normal University, Jinhua 321004, P.R. China}

\thanks[Y]{Corresponding author.}

\baselineskip 20pt

\date{}
\maketitle
 \vspace{4mm}

\begin{abstract}
Let  $X$ be a finite set and  $\mathfrak p\subseteq 2^X$, the power
set of $X$, satisfying  three conditions: (a) $\mathfrak p$ is an
ideal in $2^X$, that is, if $A\in \mathfrak p$ and $B\subset A$,
then $B\in \mathfrak p$; (b) For $A\in 2^X$ with $|A|\geq 2$, $A\in
\mathfrak p$ if $\{x,y\}\in \mathfrak p$ for any $x,y\in A$ with
$x\neq y$; (c)  $\{x\}\in \mathfrak p$ for every $x\in X$. The pair
$(X,\mathfrak p)$ is called a symmetric system if there is a group
$\Gamma$ transitively acting on $X$ and preserving the ideal
$\mathfrak p$.  A family $\{A_1,A_2,\ldots,A_m\}\subseteq 2^X$ is
said to be a cross-$\mathfrak{p}$-family of $X$ if  $\{a, b\}\in
\mathfrak{p}$ for any $a\in A_i$ and $b\in A_j$ with $i\neq j$. We
prove that if $(X,\mathfrak p)$ is a symmetric system and
$\{A_1,A_2,\ldots,A_m\}\subseteq 2^X$ is a
cross-$\mathfrak{p}$-family of $X$, then
\[\sum_{i=1}^m|{A}_i|\leq\left\{
\begin{array}{cl}
|X| & \hbox{if $m\leq \frac{|X|}{\alpha(X,\, \mathfrak p)}$,} \\
m\, \alpha(X,\, \mathfrak p) & \hbox{if $m\geq
\frac{|X|}{\alpha{(X,\, \mathfrak p)}}$,}
\end{array}\right.\]
where $\alpha(X,\, \mathfrak p)=\max\{|A|:A\in\mathfrak p\}$. This
generalizes Hilton's theorem on cross-intersecting families of
finite sets, and provides analogs for cross-$t$-intersecting
families of finite sets, finite vector spaces and permutations, etc.
 Moreover, the primitivity of symmetric systems is introduced to
characterize the optimal families.

\begin{keyword} intersecting family, cross-intersecting family,
symmetric system, Erd\H{o}s-Ko-Rado theorem
 \\[7pt]
\end{keyword}
\end{abstract}
\end{frontmatter}

\newcommand{\lr}[1]{\langle #1\rangle}
\newcommand{\qchoose}[2]{{ #1   \atopwithdelims[]  #2 }}

\parindent 17pt
\baselineskip 17pt
\section{Introduction}

A family $\mathcal A$ of sets is said to be intersecting if $A\cap
B\neq \emptyset$  for any $ A, B\in\mathcal  A$. A classical result
on intersecting families is due to Erd\H{o}s, Ko and Rado, which
says that if  $\mathcal A$ is an intersecting family consisting of
$k$-element subsets of an $n$-element set with $n\geq 2k$, then
$|\mathcal A|\leq {n-1\choose k-1}$, and if $n>2k$, equality
holds if and only if every subset in $\mathcal A$ contains a fixed
element.

The Erd\H{o}s-Ko-Rado theorem has many generalizations, analogs and
variations. First, the notion of intersection is generalized to
$t$-intersection, and finite sets  are analogous to finite vector
spaces, permutations and other mathematical objects. Second, intersecting families are
generalized to cross-intersecting families: $\mathcal A_1,\mathcal
A_2, \ldots, \mathcal A_m$ are said to be cross-intersecting if
$A\cap B\neq\emptyset$ for any $A\in\mathcal A_i$ and $B\in\mathcal
A_j$, $i\neq j$. Clearly, if $\mathcal A_1=\mathcal A_2= \ldots=
\mathcal A_m=\mathcal A$, then $\mathcal A$ is an intersecting family. Combining the two points of view, we may
consider the cross-$t$-intersecting families over finite vector
spaces,  permutations, etc.

A nice result on cross-intersecting families is given by Hilton
\cite{hilton} as follows.
\begin{theorem}(Hilton \cite{hilton})\label{hil}
Let $\mathcal{A}_1,\mathcal{A}_2,\ldots,\mathcal{A}_m$ be
cross-intersecting families of $k$-element subsets of an $n$-element
set $X$ with $\mathcal{A}_1\neq\emptyset$. If  $k\leq n/2$, then
\begin{eqnarray}\label{hilton}
 \sum_{i=1}^m|\mathcal{A}_i|\leq\left\{
                                        \begin{array}{cl}
                                          \binom{n}{k}, & \hbox{if $m\leq \frac{n}{k}$;} \\
                                          m\binom{n-1}{k-1}, & \hbox{if $m\geq \frac{n}{k}$.}
                                        \end{array}
                                      \right.
\end{eqnarray}
Unless $m=2=n/k$, the bound is attained if and only if one of the
following holds:
\begin{enumerate}
  \item [\rm (i)] $m<n/k$ and
$\mathcal{A}_1=\{A\subset X: |A|=k\}$,  and
$\mathcal{A}_2=\cdots=\mathcal{A}_m=\emptyset$;
  \item [\rm (ii)] $m>n/k$ and
$|\mathcal{A}_1|=|\mathcal{A}_2|=\ldots=|\mathcal{A}_m|=\binom{n-1}{k-1}$;
  \item [\rm (iii)] $m=n/k$ and
$\mathcal{A}_1,\mathcal{A}_2,\ldots,\mathcal{A}_m$ are as in \rm(i)
or (ii).
\end{enumerate}
\end{theorem}

Recently, Borg gives a simple proof of the above theorem
\cite{pborg1}, and generalizes it to   labeled sets \cite{borg}
and  permutations \cite{pborg}. Inspired by his proofs we shall
present a general result on  cross-intersecting, or
cross-$t$-intersecting families of finite sets, finite vector spaces, permutations, etc.  To do this, we introduce a general
definition.

Let  $X$ be a finite set and  $\mathfrak p\subseteq 2^X$, the power
set of $X$, satisfying  three conditions as follows:
\begin{enumerate}
  \item[]
\begin{itemize}
  \item [(a)]$\mathfrak p$ is an ideal in $2^X$, that is, if $A\in \mathfrak p$ and $B\subset A$, then $B\in \mathfrak p$;
  \item [(b)] For $A\in 2^X$ with $|A|\geq 2$, $A\in \mathfrak p$ if  $\{x,y\}\in \mathfrak p$ for any $x,y\in A$ with $x\neq y$;
\item [(c)]  $\{x\}\in \mathfrak p$ for every $x\in X$.
\end{itemize}
\end{enumerate}

 Note that condition (a) is essential and (c) is to avoid
  trivial cases. If
ignore conditions (b) and (c), the pair $(X,\mathfrak p)$ is  an
(abstract) simplicial complex in topology, or a hereditary family in
extremal set theory (see e.g. \cite[p.86]{Engel97} or \cite{borg5}).
If ignore (b), $\mathfrak p$ is called a full hereditary family in
 \cite[p.86]{Engel97}.
   Condition (b) is not redundant in most discussions on extremal
combinatorics, and is necessary in our argument.

 Clearly, $\mathfrak p$ defines a binary relation
``$\sim_\mathfrak p$" on $X$: $x\sim_\mathfrak p y$ if and only if
$\{x,y\}\in \mathfrak p$ for any $x,y\in X$. This relation is
reflexive and symmetric, i.e., $x\sim_\mathfrak p x$ for every $x\in
X$, and  $x\sim_\mathfrak p y$ implies $y\sim_\mathfrak p x$.
Conversely, given a reflexive and symmetric  binary relation
``$\sim$" on $X$, we can get an ideal $\mathfrak p$ in $2^X$:
$A\subset X$ is in $\mathfrak p$ if $a\sim b$ for any $a,b\in A$.
Moreover, $\mathfrak p$ also defines a property on $2^X$: a subset
$A$ of $X$ has the property $\mathfrak p$ if $A\in \mathfrak p$.
Therefore, we call the pair $(X,\mathfrak p)$ a {\it $\mathfrak
p$-system}, or a system, for short.

 An element of $\mathfrak{p}$ is
also called a $\mathfrak{p}$-subset of $X$. A family
$\{A_1,A_2,\ldots,A_m\}\subseteq 2^X$ is said to be a {\it
cross-$\mathfrak{p}$-family of $X$} if $\{a, b\}\in \mathfrak{p}$
for any $a\in A_i$ and $b\in A_j$ with $i\neq j$. By definition we
see that if $\{A_1,A_2,\ldots,A_m\}$ is a
cross-$\mathfrak{p}$-family and $A_1=A_2=\cdots=A_m=A$, then $A$ is
a $\mathfrak{p}$-subset. Write
$$\alpha(X,\mathfrak
p):=\max\{|A|:A\in\mathfrak p\}$$
 and  $$\alpha_m(X,\mathfrak
p):=\max\left\{\sum_{i=1}^m|A_i|: \{A_1,A_2,\ldots,A_m\}\ \mbox{is a
cross-$\mathfrak p$-family}\right\}.$$ A cross-$\mathfrak p$-family
$\{A_1,A_2,\ldots,A_m\}$ is said to be {\it optimal}\ if
$\sum_{i=1}^m|A_i|=\alpha_m(X,\mathfrak p)$.

We call a system  $(X,\mathfrak p)$ {\it symmetric} if there is a
group $\Gamma$ transitively acting on $X$ and preserving the
property $\mathfrak p$, i.e., for every pair $a,b\in X$ there is a
$\gamma\in \Gamma$ such that $b=\gamma(a)$, and $A\in \mathfrak p$
implies $\delta(A)\in \mathfrak p$ for every $\delta\in \Gamma$. In
this case we say that the group $\Gamma$ transitively acts on
$(X,\mathfrak p)$.

 Two typical examples of symmetric systems are as follows.

\begin{example}\label{subset}\rm For a positive integer $n$, let $[n]$ denote the set
$\{1,2,\ldots,n\}$. By $\mathcal C_{n}^{k}$ we denote the set of all
$k$-element subsets of $[n]$, as known for $\binom{[n]}{k}$  in many
literatures.  Then $|\mathcal C_{n}^{k}|={n\choose k}$. A subset
$\mathcal A$ of $\mathcal C_{n}^{k}$ is said to be a
$t$-intersecting family if $|A\cap B|\geq t$ for any $A,B\in
\mathcal A$, where $1\leq t\leq k$. For convenience, we regard the
empty set  as a $t$-intersecting family. Let $\mathfrak i_t$ be the
collection of all $t$-intersecting families in $\mathcal C_{n}^{k}$.
Then, it is clear that $\mathfrak i_t$ is an ideal of the power set
of  $\mathcal C_{n}^{k}$, and satisfies condition (b). When $t=1$,
$\mathfrak i_t$ is abbreviated as $\mathfrak i$. The
Erd\H{o}s-Ko-Rado theorem and Theorem \ref{hil} say that
$\alpha(\mathcal C_{n}^{k}, \mathfrak i)={n-1\choose k-1}$  and
$\alpha_m(\mathcal C_{n}^{k}, \mathfrak i)=\max\left\{{n\choose k},
m{n-1\choose k-1}\right\}$ for $n\geq 2k$, respectively. In fact,
Erd\H{o}s, Ko and Rado \cite{EKR} also proved $\alpha(\mathcal
C_{n}^{k}, \mathfrak i_t)=\binom{n-t}{k-t}$ for $t>1$ and $n\geq
n_0(k,t)$, a sufficiently large positive integer depending on $k$
and $t$. The smallest $n_0(k,t)= (k-t+1)(t+1)$ was determined by
Frankl \cite{pf} for $t\geq 15$ and subsequently determined by
Wilson \cite{wilton} for all $t$. It is well known that the
symmetric group $S_n$ transitively acts on $\mathcal C_{n}^{k}$ in a
natural way, and preserves $\mathfrak i_t$. Therefore, $(\mathcal
C_{n}^{k},\mathfrak i_t)$ is symmetric.
\end{example}
\begin{example}\label{subspace}\rm Let $\mathcal L_{n,k}(q)$ denote the set of all $k$-dimensional
subspaces of an $n$-dimensional vector space over a $q$-element
field. Then  $|\mathcal L_{n,k}(q)|=\qchoose
nk=\frac{\{n\}!}{\{k\}!\{n-k\}!}$ where $\{k\}=1+q+\cdots+q^{k-1}$
and $\{k\}!=\{k\}\{k-1\}\cdots\{1\}$.
 A subset
$\mathcal A$ of $\mathcal L_{n,k}(q)$ is said to be a
$t$-intersecting family if $\dim(A\cap B)\geq t$ for any $A,B\in
\mathcal A$, where $1\leq t\leq k$. We still use $\mathfrak i_t$ to
denote the collection of all $t$-intersecting families in $\mathcal
L_{n,k}(q)$, and abbreviate $\mathfrak i_1$ as $\mathfrak i$. That
$\alpha(\mathcal L_{n,k}(q),\mathfrak i)=\qchoose{n-1}{k-1}$ was
first established by Hsieh \cite{Hsieh} for  $k<n/2$, and by Greene
and  Kleitman \cite{Greene-Kleit} for $k|n$. For $t\geq 2$, Frankl
and Wilson \cite{Frankl-Wilson} proved that $\alpha(\mathcal
L_{n,k}(q),\mathfrak i_t)=\max\left\{\qchoose
{n-t}{k-t},\qchoose{2k-t}{k}\right\}$ for $n\geq 2k-t$. Analogously
to  $(\mathcal C_{n}^{k},\mathfrak i_t)$, the general linear group
$GL(n,q)$ transitively acts on $\mathcal L_{n,k}(q)$ and preserves
$\mathfrak i_t$. Therefore, $(\mathcal L_{n,k}(q),\mathfrak i_t)$ is
also symmetric.
\end{example}

To our knowledge, there is no information on  $\alpha_m(\mathcal
C_{n}^{k}, \mathfrak i_t)$ for $t>1$ and $\alpha_m(\mathcal
L_{n,k}(q),\mathfrak i_t)$ for $t\geq 1$.

In this paper we shall generalize  Theorem \ref{hil}  to all
symmetric systems $(X,\mathfrak p)$ up to $\alpha(X,\mathfrak{p})$.
The main result will be presented in the next section. To
characterize the optimal cross-$\mathfrak p$-families we introduce
the primitivity of the symmetric systems, and give its main
characters in Section 3. As applications of results in Section 3, we
prove in Section 4 that the symmetric systems defined on finite
sets, finite vector spaces and symmetric groups are all primitive
except a few trivial cases.

\section{Cross-intersecting families of  symmetric systems}

Given a system $(X,\mathfrak p)$, we can construct a simple graph,
written as $G(X,\mathfrak p)$, whose vertex set is $X$, and
$\{a,b\}$ is an edge if $\{a,b\}\not\in \mathfrak p$. Then every
subset of $X$ in $\mathfrak p$ corresponds to an independent set of
$G(X,\mathfrak p)$. Conversely, given a simple graph $G$, we obtain
a system $(X(G),\mathfrak p(G))$, where $X(G)$ is the vertex set
$V(G)$ of $G$ and $\mathfrak p(G)$ consists of all independent sets
of $G$. It is clear that $\alpha(X(G),\mathfrak p(G))=\alpha(G)$,
the independence number of $G$.

By $I(X, \mathfrak p)$ we denote the set of all maximal-sized
$\mathfrak p$-subsets of $X$. Similarly, for a graph $G$, let $I(G)$
 denote the set of all  maximal-sized independent sets of $G$. For $B\subseteq V(G)$, let $G[B]$
denote the induced subgraph of $G$ by $B$.

 The notations introduced below have graph-theoretic intuition.

Let $(X, \mathfrak p)$ be a $\mathfrak p$-system.   For $B\subseteq
X$, we abbreviate $\alpha(B, \mathfrak p\cap 2^B)$ as $\alpha(B,
\mathfrak p)$. Clearly, $\alpha(B, \mathfrak p)$ equals
 $\alpha(G[B])$, where $G=G(X,\mathfrak p)$.
For $A\subseteq X$, set
$$
N_{X, \mathfrak p}[A]=A\cup\{b\in X:\mbox{ $\{a,b\}\not\in \mathfrak
p$ for some $a\in A$ }\}$$ and
$$  \bar N_{X,
\mathfrak p}[A]=X-N_{X, \mathfrak p}[A].$$ If there is no
possibility of confusion, we abbreviate $N_{X, \mathfrak p}[A]$ as
$N[A]$. From definition we see that $N[\emptyset]=\emptyset$; $
N[A]=X$ if $A\in I(X, \mathfrak p)$; if both $B\subseteq A$ and
$C\subseteq \bar N[A]$ are in $\mathfrak p$, then $B\cup C\in
\mathfrak p$.

  We call $(X, \mathfrak p)$
\textit{connected (disconnected)} if the graph $G(X, \mathfrak p)$
is connected (disconnected). By definition we see that $(X,
\mathfrak p)$ is disconnected if and only if there is a proper
subset $A\subset X$ such that  $\bar N[A]=X-A$, and,  $(X,\mathfrak
p)$ is symmetric if and only if $G(X,\mathfrak p)$ is
vertex-transitive.

In the context of vertex-transitive graphs, the ``No- Homomorphism"
lemma  is useful to get bounds on the size of independent sets.
\begin{lemma}\label{ac}(
Albertson and Collins \cite{makl}) Let $G$ and $H$ be two graphs
such that $G$ is vertex-transitive and there exists a homomorphism
$\phi: H\mapsto G$. Then $\frac{\alpha(G)}{|V(G)|}\leq
\frac{\alpha(H)}{|V(H)|}$, and equality holds if and only if  for
each $I\in I(G)$, $\phi^{-1}(I)\in I(H)$.
\end{lemma}

In the above lemma, by taking $H$  as an induced subgraph of $G$ and
$\phi$ as the embedding mapping, we obtain the following theorem,
which is more convenient in our argument.

\begin{theorem}\label{caku}(Cameron and Ku \cite{Cameron})
Let $G$ be a vertex-transitive graph and  $B$  a subset of $V(G)$.
Then any independent set  $S$ in $G$ satisfies that
$\frac{|S|}{|V(G)|}\leq \frac{\alpha(G[B])}{|B|}$, equality implies
that $|S\cap B|=\alpha(G[B])$.
\end{theorem}

In \cite{zhhj}, the second author of this paper proved Lemma
\ref{cl1} and Theorem \ref{irregular} below in terms of graph
theory. He also introduced the concept of imprimitive independent
sets of a vertex-transitive graph.  For completeness we restate them
in terms of symmetric systems and provide proofs for them.

\begin{lemma}\label{cl1} Let $(X, \mathfrak p)$ be a symmetric
system. Then $\frac{|A|}{|N[A]|}\leq \frac{\alpha(X, \mathfrak
p)}{|X|}$ for an arbitrary $\mathfrak p$-subset $A$ of $X$. Equality
implies that $|S\cap N[A]|=|A|$  for every $S\in
  I(X, \mathfrak p)$,
  and $\frac{\alpha(\bar N[A]\!,\, \mathfrak p)}{|\bar N[A]|}=\frac{\alpha(X\!,\, \mathfrak p)}{|X|}$.
\end{lemma}
\textbf{Proof.}  Let $C$ be a maximal-sized $\mathfrak p$-subset of
$\bar N[A]$. Clearly, $A\cup C$ is a $\mathfrak p$-subset of $X$ and
$$\frac{|A\cup C|}{|X|}=
\frac{|A|+\alpha(\bar N[A], \mathfrak p)}{|N[A]|+|\bar N[A]|}\leq
\frac{\alpha(X,\mathfrak p)}{|X|}.$$
 Since $\frac{\alpha(\bar N[A],\, \mathfrak p)}{|\bar N[A]|}\geq
\frac{\alpha(X,\, \mathfrak p)}{|X|}$  by  Theorem  \ref{caku},
$\frac{|A|}{|N[A]|}\leq \frac{\alpha(X,\, \mathfrak p)}{|X|}$.
Equality implies that  $\frac{\alpha(\bar N[A],\, \mathfrak
p)}{|\bar N[A]|}=\frac{\alpha(X,\, \mathfrak p)}{|X|}$ and
$\alpha(X, \mathfrak p)=\alpha(\bar N[A],\mathfrak p)+|A|$. Again by
Theorem \ref{caku}, we have that
 $|S\cap \bar N[A]|=|\alpha(\bar N[A], \mathfrak p)|$ and
$|S|=|S\cap N[A]|+|S\cap \bar N[A]|$ for every $S\in I(X, \mathfrak
p)$. Therefore,  $|S\cap N[A]|=|A|$ for every $S\in I(X, \mathfrak
p)$, completing the proof. \qed

 In \cite{zhhj}, a graph $G$ is called {\em IS-imprimitive} (independent-set-imprimitive) if
  there is an independent set $A$  of  $G$  such that $|A|< \alpha(G)$ and
$\frac{|A|}{|N[A]|}=\frac{\alpha(G)}{|V(G)|}$, and $A$ is called an
{\em imprimitive independent set} of $G$. In any other case, $G$ is
called {\em IS-primitive}. In this paper, we say a system
$(X,\mathfrak p)$ is {\em $\mathfrak p$-imprimitive ($\mathfrak
p$-primitive)} if the graph $G(X,\mathfrak p)$ is IS-imprimitive
(IS-primitive); a $\mathfrak p$-subset $A$ is called imprimitive if
$A$ is an imprimitive independent set of $G(X,\mathfrak p)$.
 From definition we
 see that a disconnected symmetric system $(X,\mathfrak p)$ is $\mathfrak p$-imprimitive and hence a ${\mathfrak p}$-primitive symmetric system $(X,\mathfrak p)$ is connected.

We now  contribute to $\alpha_m(X,\mathfrak p)$. Note that in a
series of papers \cite{borg,pborg1,pborg,borg6} Borg determined this
value for various cross-intersecting families. An important step in
his proofs was inequality (\ref{trik}) below he established for some
special intersecting families. We find that the inequality for
$\mathfrak p$-subsets in  symmetric systems is a consequence of
Theorem \ref{caku}, stated as follows.

\begin{corollary}\label{ab}
Let $(X, \mathfrak p)$ be a symmetric system, and let $A$ be a
$\mathfrak p$-subset of $X$. Then
\begin{equation}\label{trik}
    |A|+\frac{\alpha(X, \mathfrak p)}{|X|}|\bar{N}[A]|\leq \alpha(X,
\mathfrak p).
\end{equation}
Equality holds if and only if  $A=\emptyset$ or $|A|=\alpha(X,
\mathfrak p)$ or
 $A$ is  an imprimitive $\mathfrak p$-subset.
\end{corollary}
\textbf{Proof.} If  $A=\emptyset$ or $|A|=\alpha(X, \mathfrak p)$,
equality trivially holds.  Suppose that $0<|A|<\alpha(X, \mathfrak
p)$ and  $B$ is a maximal-sized $\mathfrak p$-subset  in
$\bar{N}[A]$, that is, $|B|=\alpha(\bar{N}[A], \mathfrak p)$. Then
$A\cup B$ is also a $\mathfrak p$-subset of $X$, so $|A|+|B|\leq
\alpha(X, \mathfrak p)$, and Theorem \ref{caku} implies that
$\frac{|B|}{|\bar{N}[A]|}\geq \frac{\alpha(X, \mathfrak p)}{|X|}$.
Therefore, $$\label{re}|A|+\frac{\alpha(X, \mathfrak
p)}{|X|}|\bar{N}[A]|\leq |A|+|B|\leq \alpha(X, \mathfrak p). $$ If
$\alpha(X, \mathfrak p)=|A|+\frac{\alpha(X, \mathfrak
p)}{|X|}|\bar{N}[A]|= |A|+\frac{\alpha(X, \mathfrak
p)}{|X|}(|X|-|N[A]|)$, then $\frac{|A|}{|N[A]|}=\frac{\alpha(X,
\mathfrak p)}{|X|}$, i.e., $A$ is an imprimitive $\mathfrak
p$-subset. \qed

The following theorem is  the main result of this paper.
\begin{theorem}\label{m}
Let $(X, \mathfrak p)$ be a connected symmetric system, and let
$\{{A}_1, {A}_2,\ldots,{A}_m\}$ be a cross-$\mathfrak p$-family over
$X$ with $A_1\neq\emptyset$. Then
$$ \sum_{i=1}^m|{A}_i|\leq\left\{
                                        \begin{array}{cl}
                                          |X| & \hbox{if $m\leq \frac{|X|}{\alpha(X,\, \mathfrak p)}$;} \\
                                          m\,\alpha(X,\,  \mathfrak p) & \hbox{if $m\geq \frac{|X|}{\alpha{(X,\,  \mathfrak p)}}$,}
                                        \end{array}
                                      \right.
                                      $$
and the bound is attained if and only if one of the following holds:
\begin{enumerate}
  \item [\rm(i)] $m<\frac{|X|}{\alpha(X,\,  \mathfrak p)}$ and $A_1=X$, $A_2=\ldots=A_m=\emptyset$,
  \item [\rm(ii)] $m>\frac{|X|}{\alpha(X,\,  \mathfrak p)}$ and ${A}_1=\ldots={A}_m=I\in I(X, \mathfrak p)$,
  \item [\rm(iii)] $m=\frac{|X|}{\alpha(X,\,  \mathfrak p)}$ and either ${A}_1,{A}_2,\ldots,{A}_m$ are as in {\rm(i) or
(ii)}, or there is an imprimitive $\mathfrak p$-subset $A$ such that
$A\subseteq {A}_i$, $i=1,2,\ldots,m$, and $\{A_1',A_2',\ldots,A_m'\}$ is a cross-$\mathfrak p$-family and a partition of
$\bar N[A]$, where
$A_i'=A_i-A$, $i=1,2\ldots,m$.
\end{enumerate}
\end{theorem}
{\bf Proof. } Following Borg's notation  in
\cite{pborg1,pborg,borg6}, write $A_i^*=\{a\in A_i: \{a,b\}\in
\mathfrak p\ \mbox{for every $b\in A_i$}\}$, $A_i'=A_i-A_i^*$,
$A^*=\displaystyle{\cup_{i=1}^{m}}{A}_i^*$ and
$A'=\displaystyle{\cup_{i=1}^{m}}{A}_i'$. It is clear that $A^*$ is
a $\mathfrak p$-subset and  ${A}'\subseteq \bar{N}[A^*]$. From
definition it follows that $A_i\cap A_j\subseteq A_i^*\cap A_j^*$,
therefore $A_i'\cap A_j'=\emptyset$ for $i\neq j$, thus
$|{A}'|=\sum_{i=1}^m|{A}'_i|$. By Corollary \ref{ab} we have that
\begin{eqnarray*}
\sum_{i=1}^m|{A}_i|&=&\sum_{i=1}^m|{A}'_i|+\sum_{i=1}^m|{A}_i^*|\leq|{A}'|+m|A^*|\leq |\bar N[A^*]|+m|A^*|\\
&=&\frac{|X|}{\alpha(X, \mathfrak p)}\left(\frac{\alpha(X, \mathfrak
p)}{|X|}|\bar N[A^*]|+|A^*|\right)
+\left(m-\frac{|X|}{\alpha(X, \mathfrak p)}\right)|A^*|\\
&\leq&|X|+\left(m-\frac{|X|}{\alpha(X, \mathfrak p)}\right)|A^*|.
\end{eqnarray*}

If $m<\frac{|X|}{\alpha(X, \mathfrak p)}$, then
$\sum_{i=1}^m|{A}_i|\leq |X|$, and equality  implies
$A^*=\emptyset$, hence ${A_i}={A_i}'$ for every $i\in[m]$, and we
thus have that the corresponding graph $G(X, \mathfrak p)$ is a
union of the induced subgraphs $G(X, \mathfrak p)[A_i']$'s. Then,
the connectivity of $(X, \mathfrak p)$ yields that one of them is
$X$ and the others are empty, as (i).

If $m>\frac{|X|}{\alpha(X, \mathfrak p)}$, then
$\sum_{i=1}^m|{A}_i|\leq m\,\alpha(X, \mathfrak p)$ and  equality
implies that ${A}_1^*=\cdots=A_m^*=A^*$ and $|A^*|=\alpha(X, \mathfrak
p)$, as (ii).

If $m=\frac{|X|}{\alpha(X, \mathfrak p)}$, then
$\sum_{i=1}^m|{A}_i|\leq |X|$, and equality implies that
${A}_1^*=\cdots=A_m^*=A^*$ and $\frac{\alpha(X, \mathfrak
p)}{|X|}|\bar N[A^*]|+|A^*|=\alpha(X, \mathfrak p)$. Then Corollary
\ref{ab} implies that $|A^*|=0$ or $|A|=\alpha(X, \mathfrak p)$ or
$A^*$ is an imprimitive $\mathfrak p$-subset. In the last case,
$\{A_1', A_2', \dots, A_m'\}$ is a cross-$\mathfrak p$-family, and a
 partition of $\bar N[A^*]$. \qed

From the above theorem we see that if $(X, \mathfrak p)$ is
symmetric and $\mathfrak p$-primitive (hence connected), then
$\alpha_m(X,\mathfrak p)$ is uniquely determined by
$\alpha(X,\mathfrak p)$, i.e.,
\[\alpha_m(X,\mathfrak p)=\max\left\{|X|,m\,\alpha(X,\mathfrak
p)\right\},\] and an optimal cross-$\mathfrak p$-family is one of
the forms $\{X,\emptyset,\ldots,\emptyset\}$ and $\{A,A,\ldots,A\}$
where $A\in \mathfrak p$ with $|A|=\alpha(X,\mathfrak p)$.

For the $(X,\mathfrak p)$ dealt with in this field, however,
$\alpha(X,\mathfrak p)$ is usually well known, and the symmetric
property of $(X,\mathfrak p)$ is easy to verify. So we concentrate
on the primitivity of symmetric systems in the next two sections.

\section{Primitivity of symmetric systems}

This concept comes from permutation groups. Let $X$ be a set, and
$\Gamma$ a group transitively acting on $X$.  Then $\Gamma$ is said
to be imprimitive  on $X$ if it preserves a nontrivial partition of
$X$, called a block system, each element of which is called a block.
In any other case $\Gamma$ is primitive on $X$. More precisely,
$\Gamma$ is imprimitive  on $X$ if  there is nontrivial partition
$X=\cup_{i=1}^kX_i$ such that  $\gamma(X_i)$ is a block of the
partition for every $\gamma\in \Gamma$ and $i=1,2,\ldots,k$. Here
$\gamma(X_i)$ denotes the set $\{\gamma(x):x\in X_i\}$.

A classical result on the primitivity of group actions is the
following theorem (cf. \cite[Theorem 1.12]{BA}).

\begin{theorem}\label{max-prim} Suppose that a group $\Gamma$  transitively acts on $X$. Then  $\Gamma$  is
primitive on $X$ if and only if for each $a\in X$, $\Gamma_a$ is a
maximal subgroup of $\Gamma$. Here $\Gamma_a=\{\gamma\in
\Gamma:\gamma(a)=a\}$,  the stabilizer of $a\in X$.
\end{theorem}

 The following theorem explains why  a symmetric system is called  primitive
or imprimitive.

\begin{theorem}\label{irregular}
Let $(X, \mathfrak p)$ be an imprimitive symmetric system, $A$ a
maximal-sized imprimitive $\mathfrak p$-subset of $X$, $D=X-N[A]$,
and let $\Gamma$ be the group transitively acting on $(X,\mathfrak
p)$. Then $\frac{\alpha(D\!,\, \mathfrak p)}{|D|}=\frac{\alpha(X\!,
\,\mathfrak p)}{|X|}$ and $\{\sigma(D):\sigma\in \Gamma\}$ forms a
partition of $X$.
\end{theorem}
\textbf{Proof.} First, suppose that $A$ and $B$ are two imprimitive
$\mathfrak p$-subsets of $X$, and write $C=A\cup (B-N[A])$. We claim
that $C$ is a $\mathfrak p$-subset satisfying $N[C]=N[A]\cup N[B]$
and $\frac{|C|}{|N[C]|}=\frac{\alpha(X, \mathfrak p)}{|X|}$.

 To prove this claim we write $ N[A]\cup N[B]=M$. From definition it is
easily seen that  $C$ is also a $\mathfrak p$-subset and
$N[C]\subseteq M$. Since $\frac{|B|}{|N[B]|}=\frac{\alpha(X,
\mathfrak p )}{|X|}$, by Lemma \ref{cl1} we have that $|S\cap
N[B]|=|B|$ for all $S\in I(X, \mathfrak p)$. So, $B\cup (S-N[B])$ is
also a maximal-sized $\mathfrak p$-subset of $X$ for every $S\in
I(X, \mathfrak p)$. By repeating this process for the maximal-sized
$\mathfrak p$-subset $B\cup (S-N[B])$  and the imprimitive
$\mathfrak p$-subset $A$ we have that
\begin{eqnarray*}
  && A\cup((B\cup (S-N[B]))-N[A])
\\
&=&A\cup (B-N[A])\cup ((S-N[B])- N[A])=C\cup (S-M)
\end{eqnarray*}
is also a maximal-sized $\mathfrak p$-subset of $X$, which implies
that $|S\cap M|=|C|$ for every $S\in I(X, \mathfrak p)$. Given a
$u\in X$, suppose there are $r$ maximal-sized $\mathfrak p$-subsets
containing $u$. Since $(X, \mathfrak p)$ is symmetric, it is easily
seen that the number $r$ is independent on the choice of $u$. Let us
count pairs $(x,S)$ with $x\in M\cap S,\ S\in I(X,\mathfrak p)$, in
two ways.  Since  $|M\cap S|=|C|$ for every $S\in I(X, \mathfrak
p)$, the number of the pairs is clearly equal to $|C||I(X, \mathfrak
p)|$. On the other hand, for each $x\in M$ there are $r$ $S$'s in
$I(X, \mathfrak p)$ with $x\in S$. So the number is also equal to
$r|M|$, proving  $r|M|=|C||I(X, \mathfrak p)|$. Similarly, by
counting pairs $(x,S)$ with $x\in S\in I(X,\mathfrak p)$ in two ways
we obtain $r|X|=\alpha(X, \mathfrak p)|I(X, \mathfrak p)|$.
Combining the above two equalities gives
$\frac{|C|}{|M|}=\frac{\alpha(X, \mathfrak p)}{|X|}$. Thus, by Lemma
\ref{cl1} we have that
$$\frac{\alpha(X, \mathfrak p)}{|X|}\geq \frac{|C|}{|N[C]|}\geq
\frac{|C|}{|M|}=\frac{\alpha(X, \mathfrak p)}{|X|}.$$ Hence $N[C]=M$
and $\frac{|C|}{|N[C]|}=\frac{\alpha(X, \mathfrak p)}{|X|}$, proving
our claim.

We now close the proof of the theorem. Let $A$ be a maximal-sized
imprimitive $\mathfrak p$-subset of $X$. From definition it follows
that $N[\sigma(A)]=\sigma(N[A])$ for all $\sigma\in \Gamma$. Suppose
that there exists a $\sigma\in \Gamma$ such that $ \sigma(D)\neq D$
and $\sigma(D)\cap D\neq \emptyset$. Then $\sigma(N[A])\neq N[A]$,
hence $|N[A]\cup \sigma\big(N[A]\big)|>|N[A]|$. Set $A'=A\cup
(\sigma(A)-N[A])$. Then $A'$ is also a $\mathfrak p$-subset of $X$.
By the above claim  we have that $N[A']=N[A]\cup
\sigma\big(N[A]\big)$ and $\frac{|A'|}{|N[A']|}=\frac{\alpha(X,
\mathfrak p)}{|X|}=\frac{|A|}{|N[A]|}$, which implies $|A'|> |A|$.
On the other hand, from definition  it follows that each element of
$\sigma(D)\cap D$ does not belong to $ N[A]\cup
\sigma\big(N[A]\big)$, so $ N[A']\neq X$, yielding $|A'|< \alpha(X,
\mathfrak p)$. It contradicts the maximality of $A$, thus proving
that  $\sigma(D)= D$ or $\sigma(D)\cap D= \emptyset$ for each
$\sigma\in \Gamma$. The transitivity of $\Gamma$ on $X$ implies that
$X=\cup_{\sigma\in\Gamma}\sigma(D)$. Furthermore,  for any
$\sigma,\gamma\in \Gamma$, if $\sigma(D)\cap\gamma(D)\neq\emptyset$,
then $(\gamma^{-1}\sigma)(D)\cap D\neq\emptyset$, implying
$(\gamma^{-1}\sigma)(D)=D$, i.e., $\sigma(D)=\gamma(D)$. Therefore,
$\{\sigma(D):\sigma\in \Gamma\}$ is a partition of $X$. \qed

By  Theorem \ref{irregular}  and Theorem \ref{max-prim} we obtain
the following consequences.

\begin{corollary}\label{coropri}
Suppose that a group $\Gamma$ transitively acts on $(X, \mathfrak
p)$. Then $(X, \mathfrak p)$ is $\mathfrak p$-primitive if one of
the following conditions holds.
\begin{enumerate}
  \item[]
  \begin{itemize}
    \item[\rm (i)] $\Gamma$ is primitive on $X$, or equivalently, $\Gamma_a$
    is a
maximal subgroup of $\Gamma$ for each $a\in X$.
    \item[\rm (ii)] $\Gamma$ is imprimitive on $X$, but each block
    $D$ satisfies $\frac{\alpha(D,\mathfrak
p)}{|D|}> \frac{\alpha(X,\mathfrak p)}{|X|}$.
  \end{itemize}
\end{enumerate}
\end{corollary}

\section{Primitivity of some classical symmetric systems}

Finite sets, finite vector spaces and permutations  are among
 the most important finite structures in combinatorics,
especially in extremal combinatorics. In what follows we prove
the primitivity of three symmetric systems defined on them.

\begin{proposition} $(\mathcal C_n^k,\mathfrak
i_t)$ is $\mathfrak i_t$-primitive for $n\geq (k-t+1)(t+1)$ unless
$n=2k\geq 4$ and $t=1$.
 \end{proposition}
{\bf Proof}. Since the case $n\leq 3$  is trivial, we assume that
$n\geq 4$. From Example \ref{subset} we know that $(\mathcal
C_n^k,\mathfrak i_t)$ is symmetric and $\alpha(\mathcal
C_n^k,\mathfrak i_t)={n-t\choose k-t}$ for $n\geq (k-t+1)(t+1)$.
Consider the action of the symmetric group $S_n$ on $\mathcal
C_n^k$. It is well known that for each $A\in \mathcal C_n^k$, the
stabilizer $S_{n,A}$ of $A$ is isomorphic to $S_k\times S_{n-k}$,
which is a maximal subgroup of $S_n$ if $n\neq 2k$ (See e.g
\cite{sym}). Therefore, $(\mathcal C_n^k,\mathfrak i_t)$ is
$\mathfrak i_t$-primitive when $n\neq 2k$.  It is easily seen that
$\{A,[2k]-A\}$ is a block in $\mathcal C_{2k}^k$ under the action of
$S_{2k}$, and every block is of this form. On the other hand,
$\frac{\alpha(\{A,\bar A\},\mathfrak i_t)}{2}=\frac
12\geq\frac{{2k-t\choose k-t}}{{2k\choose k}}=\frac{\alpha(\mathcal
C_n^k,\mathfrak i_t)}{|\mathcal C_n^k|}$ for all $1\leq t\leq k$,
and equality holds if and only if $t=1$. By Corollary \ref{coropri},
$(\mathcal C_{2k}^k,\mathfrak i_t)$ is $\mathfrak i_t$-primitive for
$t>1$. It is clear that  $(\mathcal C_{2k}^k,\mathfrak i)$ is
disconnected, hence $\mathfrak i$-imprimitive.  \qed

\begin{proposition}  $(\mathcal L_{n,k}(q),\mathfrak
i_t)$ is $\mathfrak i_t$-primitive for all $n\geq 2k-t$.
 \end{proposition}
{\bf Proof}. It is well known \cite{Asch}  that for each $A\in
\mathcal L_{n,k}(q)$, the  stabilizer of $A$  is a maximal subgroup
of $GL(n,q)$. By Corollary \ref{coropri} $(\mathcal
L_{n,k}(q),\mathfrak i_t)$ is $\mathfrak i_t$-primitive.\qed

In the foregoing two examples, the primitivity of systems follows
directly from the primitivity of groups acting on them. However,  it
is not always the case, as we shall see.

Let us consider the set $S_n$.  A subset $A$ of $S_n$ is said to be
$t$-intersecting if any two permutations in $A$ agree in at least
$t$ points, i.e. for any $\sigma,\tau\in A$, $|\{i\in
[n]:\sigma(i)=\tau(i)\}|\geq t$. We still denote this property by
$\mathfrak i_t$.  When $t=1$, Deza and Frankl \cite{DF} showed that
a $1$-intersecting subset $A\subseteq S_n$ has size at most $(n-1)!$
and conjectured that for $t$ fixed, and $n$ sufficiently large
depending on $t$, a $t$-intersecting subset $A\subseteq  S_n$ has
size at most $(n-t)!$. Cameron and Ku \cite{Cameron} proved a
$1$-intersecting subset of size  $(n-1)!$ is a coset of the
stabilizer of a point. A few alternative proofs of Cameron and Ku's
result are given in \cite{Larose}, \cite{godsil} and \cite{wz}. To
show the transitivity of $(S_n, \mathfrak i_t)$ we consider the
action of $S_n$ on itself by the multiplication  on the left. It is
evident that the action is transitive, but is far from primitive
because the stabilizer of a point is the identity.

\begin{proposition}  $(S_n,\mathfrak
i_t)$ is $\mathfrak i_t$-primitive unless $n=3$ and $t=1$.
 \end{proposition}
\textbf{Proof.} The case $n=2$ is trivial. If $n=3$, it is easy to
verify that the graph $G(S_3, \mathfrak i)$ is disconnected and
hence $\mathfrak i$-imprimitive, while  $(S_3, \mathfrak i_t)$ for
$t=2, 3$ is $\mathfrak i_t$-primitive. We now assume that  $n\geq
4$.

We first prove that $(S_n, \mathfrak i_t)$ is connected, i.e, the
corresponding graph $G(S_n, \mathfrak i_t)$ is connected. Since
$\mathfrak i_t \subseteq \mathfrak i_1$ for $t\geq 2$, it suffices
to prove that $G(S_n, \mathfrak i)$ is connected. For any pair
$\gamma,\eta\in S_n$, let $A_j=\{i\in [n]: \eta(j)\neq i\neq
\gamma(j)\}$ for $1\leq j\leq n$. Clearly, $|A_j|\geq n-2$. For
every $J\subseteq [n]$, if $|J|=2$, then $|\cup_{j\in J}A_j|\geq
|A_j|=n-2\geq 2$. Suppose that $|J|\geq 3$. Then, for each $k\in
[n]$, since there are at most two points $i_1,i_2\in[n]$ such that
$\gamma(i_1)=\eta(i_2)=k$,  we can find a  $j\in J$ such that $k\in
A_j$, so $\cup_{j\in J}A_j=[n]$. Therefore $|\cup_{j\in J}A_j|\geq
|J|$ for all $J\subseteq [n]$. By the well-known Hall theorem
\cite{Hall} on distinct representatives of subsets, there is a
system of distinct representatives $i_1,i_2,\ldots,i_n$ for
$A_1,A_2,\ldots,A_n$. Define a permutation $\tau$ by $\tau(j)=i_j$
for $1\leq j\leq n$. It is clear that both $\{\eta,\tau\}$ and
$\{\tau,\gamma\}$ belong to $ E(G(S_n, \mathfrak i))$, proving that
$G(S_n, \mathfrak i)$ is connected.

Suppose that $(S_n, \mathfrak i_t)$ is $\mathfrak i_t$-imprimitive
for some $n\geq 4$ and $t\geq 1$. Let $A$ be a maximal-sized
imprimitive $\mathfrak i_t$-subset of $S_n$, and $D=\bar
N[A]=S_n-N[A]$. From Theorem \ref{irregular}, it follows that
$\frac{\alpha(D, \mathfrak i_t)}{|D|}=\frac{\alpha(S_n, \mathfrak
i_t)}{|S_n|}$, and $\tau D\cap D=\emptyset$ or $D$ for all $\tau\in
S_n$, and Theorem \ref{caku} implies that $|S\cap D|=\alpha(D,
\mathfrak i_t)$ for every $S\in I(S_n, \mathfrak i_t)$. Let $\sigma$
be a fixed $n$-cycle permutation in $S_n$, and  $H=\{
\sigma,\sigma^2,\ldots,\sigma^{n}=1\}$, the cyclic group generated
by $\sigma$. Then any two distinct elements of a right coset of $H$
disagree at every point. Therefore $H\rho\subset N[\{\rho\}]$ for
every $\rho\in S_n$, so $HA\subseteq N[A]$. Set $B=\{\rho \in S_n:
H\rho \subset D\}$ and  $C=\{\rho\in S_n: \mbox{$H\rho\cap N[A]\neq
\emptyset$ and $H\rho \cap D\neq \emptyset$}\}$. We now  complete
the proof by two cases.

 Case 1: $t\geq 2$. For any $\tau,\rho\in S_n$, set $F_i=F_i(\tau,\rho)=
 \{j: \tau(j)=\sigma^i\rho(j)\}$, $i=1,2,\ldots,n$.
It is easily seen that for every $j\in [n]$ there is a unique $i\in
[n]$ such that $j\in F_i$, which yields $\sum_{i=1}^n|F_i|=n$. From
this we see that there are at least half\ $F_i$'s with at most one
point, meaning that there are at least $\lceil n/2\rceil$ $i$'s such
that $\tau$ and $\sigma^i\rho$ do not agree on $t$ points. In other
words, $|H\rho \cap N[\{\tau\}]|\geq \lceil\frac n2\rceil\geq 2$,
which implies that $B=\emptyset$ and $D\subset \cup_{\rho \in
C}H\rho$. If $\sigma D\cap D\neq \emptyset$, then $\sigma D=D$,
hence $HD=D$, contradicting $B=\emptyset$. We therefore obtain that
$\sigma D\cap D=\emptyset$. Moreover, since $\frac{\alpha(\sigma D,
\mathfrak i_t)}{|\sigma D|}=\frac{\alpha( D, \mathfrak
i_t)}{|D|}=\frac{\alpha(S_n, \mathfrak i_t)}{|S_n|}$, from Theorem
\ref{caku} it follows that $|S\cap \sigma D|=\alpha(\sigma D,
\mathfrak i_t)=\alpha(D, \mathfrak i_t)$ for every $S\in I(S_n,
\mathfrak i_t)$. Note that for each $S_D\in I(D,\mathfrak i_t)$, we
have  $A\cup S_D\in I(S_n, \mathfrak i_t)$, so  $|(A\cup S_D)\cap
\sigma D|=\alpha(D, \mathfrak i_t)$. Recalling that $HA\subseteq
N[A]$, we have $$(A\cup S_D)\cap \sigma D=A\cap \sigma D\subseteq
HA\cap \sigma D=\sigma(HA\cap D)\subseteq \sigma(N[A]\cap
D)=\emptyset,$$ yielding a contradiction. Thus $(S_n, \mathfrak
i_t)$ is $\mathfrak i_t$-primitive for $t\geq 2$.

Case 2: $t=1$. By definition we see that $|A\cap H|\leq 1$. On the
other hand, from  $HA\subseteq N[A]$ and
$\frac{|A|}{|N[A]|}=\frac{\alpha(S_n, \mathfrak
i)}{|S_n|}=\frac{1}{n}$ it follows that $N[A]=HA$, that is, $N[A]$
is a union of some right cosets of $H$, so $D$ is a union of other
right cosets of $H$, i.e., $D=HB$.  By definition we also have that
$A\subseteq \bar N[D]\subseteq \bar N[H\rho]$ for every $\rho\in B$.
However, if $\tau\in \bar N[H\rho]$, i.e. $F_i(\tau,\rho)=\{j:
\tau(j)=\sigma^i\rho(j)\}\neq\emptyset$ for every $i\in [n]$,
 then $$F_i(\sigma^k\tau,\rho\}=\{j:
\sigma^k\tau(j)=\sigma^{i}\rho(j)\}=\{j:
\tau(j)=\sigma^{i-k}\rho(j)\}=F_{i-k}(\tau,\rho)\neq \emptyset$$ for
all $i,k\in[n]$ (here $i-k$ is taken to be the least positive
residue modulo $n$), therefore $H\tau\subseteq \bar N[H\rho]$. From
this it follows that $N[A]=HA\subseteq\bigcap_{\rho\in B}\bar
N[H\rho]=\bar N[D]$, which implies that  $(S_n, \mathfrak i)$ is
disconnected, yielding a contradiction. Thus $(S_n, \mathfrak i)$ is
$\mathfrak i$-primitive for $n\geq 4$.\qed

Analogously, we may consider the primitivity of symmetric systems
defined on labeled sets \cite{borg} (or signed sets \cite{borg-dm},
colored sets \cite{lw} etc) and some other permutations (see
\cite{Ku}, \cite{KuTom} and \cite{wz}).

\section*{Acknowledgements}
The authors are greatly indebted to the anonymous referees for
giving   useful comments and suggestions that have considerably
improved the manuscript. This work was partially supported by the
National Natural Science Foundation of China  (No. 10826084 and
No.10731040), Ph.D. Programs Foundation of Ministry of Education of
China (No. 20093127110001), and Innovation Program of Shanghai
Municipal Education Commission (No. 09zz134).

\end{document}